\documentclass[12pt]{article}
\usepackage{amssymb}
\usepackage{amsmath}
\usepackage{rotating}

\usepackage[a4paper]{geometry}

\newcounter{mycounter}

\newcommand{\asss}[2]{\noindent{\textbf{Assumption #1} #2}}
\newcommand{\proofs}[2]{\noindent{\textbf{Proof of #1.} #2 \hfill $\square$}}
\newtheorem{theorem}{Theorem}
\newtheorem{prop} [theorem] {Proposition}
\newtheorem{cor} [theorem] {Corollary} 

\def\plim{\mathop{\rm plim}}
\def\Pr{\mathop{\rm P}\nolimits}
\def\weakc{\implies}

\def\setitemindent#1{\settowidth{\labelwidth}{#1}%
        \leftmargini\labelwidth
        \advance\leftmargini\labelsep
   \def\@listi{\leftmargin\leftmargini
        \labelwidth\leftmargini\advance\labelwidth by -\labelsep}}

\setitemindent{(K2.4'')}

\begin{document}

\title{Model Adequacy Checks for Discrete Choice Dynamic Models\thanks{
We thank Juan Mora for useful comments. Financial support from the Fundaci\'{o}n Ram\'{o}n Areces and from the Spain Plan Nacional de I+D+I (SEJ2007-62908) is gratefully acknowledged.
}}
\author{\textsc{Igor Kheifets}\thanks{%
New Economic School, Moscow.} \textsc{\ \ and } \textsc{Carlos Velasco}\thanks{%
Department of Economics,
Universidad Carlos III de Madrid. }}
\date{February 15, 2012}
\maketitle

\abstract{
This paper proposes new parametric model adequacy tests for possibly
nonlinear and nonstationary time series models with noncontinuous data
distribution, which is often the case in applied work. In particular, we
consider the correct specification of parametric conditional
distributions in dynamic discrete choice models, not only of some particular
conditional characteristics such as moments or symmetry. Knowing the true
distribution is important in many circumstances, in particular to apply
efficient maximum likelihood methods, obtain consistent estimates of partial
effects and appropriate predictions of the probability of future events. We
propose a transformation of data which under the true conditional
distribution leads to continuous uniform iid series. The uniformity and
serial independence of the new series is then examined simultaneously. The transformation can be
considered as an extension of the integral transform tool for noncontinuous
data. We derive asymptotic properties of such tests taking into account
the parameter estimation effect. Since transformed series are iid we do not
require any mixing conditions and asymptotic results illustrate the double
simultaneous checking nature of our test. The test statistics converges
under the null with a parametric rate to the asymptotic distribution, which is
case dependent, hence we justify a parametric bootstrap approximation. The
test has power against local alternatives and is consistent. The performance
of the new tests is compared with classical specification checks for 
discrete choice models.}



\newpage

\section{Introduction}

Dynamic choice models are important econometric tools in applied
macroeconomics and finance. These are used to describe the monetary policy
decisions of central banks (Hamilton and Jord\'{a}, 2002; Basu and de Jong,
2007), for recession forecasting (Kauppi and Saikkonen, 2008; Startz, 2008)
and to model the behavior of agents in financial markets (Rydberg and
Shephard, 2003). In the simplest framework, a binary dynamic model explains
the value of an indicator variable in period $t,$ $Y_{t}\in \left\{
0,1\right\} ,$ in terms of an information set $\Omega _{t}$ available at
this period. Then $Y_{t}$ conditional on $\Omega _{t}$ is distributed as a
Bernoulli variable with expectation $p_{t}=E\left( Y_{t}|\Omega _{t}\right)
=\Pr \left( Y_{t}=1|\Omega _{t}\right) =F\left( \pi _{t}\right) $ where $\pi
_{t}=\pi \left( \Omega _{t}\right) $ summarizes the relevant information and
$F$ is a cumulative probability distribution function (cdf) monotone
increasing.
Typical specifications of the link
function $F$ are the standard normal cdf, $\Phi ,$ and the logistic cdf.

We can describe the observed values of $Y_{t}$ as $Y_{t}=1\left\{
Y_{t}^{\ast }>0\right\} $ where $Y_t^\ast$ is given by the latent variable
model
\begin{equation*}
Y_{t}^{\ast }=\pi _{t}+\varepsilon _{t}
\end{equation*}%
and $\varepsilon _{t}\sim F=F_{\varepsilon }$ are iid observations with zero
mean.

In a general specification $\pi _{t}$ is a linear combination of a set of
exogenous variables $X_{t}$ observable in $t,$ but not necessarily
contemporaneous, plus lags of $Y_{t}$ and $\pi _{t}$ itself,%
\begin{equation*}
\pi _{t}=\alpha _{0}+\alpha \left( L\right) \pi _{t}+\delta \left( L\right)
Y_{t}+X_{t}^{\prime }\beta,
\end{equation*}%
where $\delta \left( L\right) =\delta _{1}L+\cdots +\delta _{q}L^{q}$ and $%
\alpha \left( L\right) =\alpha _{1}L+\cdots +\alpha _{p}L^p.$ When $q=0,$ $%
p=1$ and $F_{\varepsilon }=\Phi $ this leads to the dynamic probit model of
Dueker (1997),%
\begin{equation*}
\pi _{t}=\pi _{0}+\delta _{1}Y_{t-1}+X_{t}^{\prime }\beta ,
\end{equation*}%
and if the roots of $1-\alpha \left( L\right) $ are out of the unit circle, $%
\pi _{t}$ can be represented in terms of infinite lags of $Y_{t}$ and $%
X_{t}. $

Many nonlinear extensions have been considered in the literature, such as
interactions with lags of $Y_{t},$ to describe the state of the economy
in the past,
\begin{equation*}
\pi _{t}=\pi _{0}+\delta _{1}Y_{t-1}+X_{t}^{\prime }\beta +\left(
Y_{t-1}X_{t}\right) ^{\prime }\gamma
\end{equation*}%
or with the sign of other variables in $X_{t},$ both stressing different
reaction functions in several regimes defined in terms of exogenous
variables at period $t$. Other specifications consider heteroskedasticity
corrections, so that Var$\left( \varepsilon _{t}\right) =\sigma ^{2}\left(
\Omega _{t}\right) ,$ for example a two regimes conditional variance, Var$%
\left( \varepsilon _{t}\right) =\sigma ^{2}\left( Y_{t-1}\right) .$

In the general ordered discrete choice model, the dependent variable takes $%
J+1$ values in a set $\mathcal{J}$, and the parametric distribution $%
\Pr (Y_{t}=j|\Omega _{t})$ can be modeled using the unobserved latent
continuous dependent variable $Y_{t}^{\ast }.$ In the typical case
where $Y_{t}=j$ if $\mu _{j-1}\leq Y_{t}^{\ast }\leq \mu _{j}$ for $j\in
\mathcal{J},$ $\mathcal{J}=\left\{ 0,1,2,...,J\right\} $ and $\varepsilon
_{t}\sim F_{\varepsilon },$ with $\mu _{-1}=-\infty $ and $\mu _{J}=\infty ,$
we have that
\begin{equation*}
\Pr (Y_{t}=j|\Omega _{t})=F_{\varepsilon }(\mu _{j}-\pi _{t})-F_{\varepsilon
}(\mu _{j-1}-\pi _{t})
\end{equation*}%
with $\alpha _{0}=0.$

Forecasting is one of the main uses of discrete choice models. In that case for the
calculation of predictions it might be necessary to resource to recursive methods when
$\delta \left( L\right) \neq 0.$ However in almost all situations parameters are unknown,
but conditional maximum likelihood (ML) estimation is straightforward given the binomial
or discrete nature of data, with typically well behaved likelihoods and asymptotic normal
estimates if the model is properly specified. The existence, representation and
probability properties of these models have been studied under general conditions by de
Jong and Woutersen (2011), who also report the consistency and asymptotic normality of ML
estimates when the parametric model is correct. However, if not, estimates will be
inconsistent and predictions can be severely biased.

This leads to the need of diagnostic and goodness-of-fit techniques, which should account for
the main features of these models, discrete nature and dynamic evolution. The first
property entails nonlinear modeling and renders invalid many methods specifically
tailored for continuous distributions. Though the latent disturbance $\varepsilon _{t}$
is continuous and with a well specified distribution, it is unobservable. Simulation
methods could be used to estimate the distribution of such innovations, but we follow an
alternative route by "\textit{continuing}" the discrete observations $Y_{t},$ so that
they have a continuous and strictly increasing conditional distribution in $\left[ -1,J\right] $ given
$\Omega _{t}$. This distribution inherits the dependence on a set of parameters and on a
conditional information set and can serve as a main tool to evaluate the appropriateness of the
hypothesized model.

Conditional distribution specification tests are often based on comparing
parametric and nonparametric estimation as Andrews' (1997) conditional
Kolmogorov test, or on the integral transform (see Bai 2003, Corradi and
Swanson 2006). The former approach is developed for different data types,
while the latter can be used only for data with  continuous
distribution. The integral transform does not require strong conditions on
the data dependence structure, so it is very useful in testing dynamic
models. However, applying the integral transform to noncontinuous data will
not bring to uniform on $[0,1]$ series, and therefore this approach can not be applied
directly to dynamic discrete choice models. To guarantee that adequacy tests
based on the integral transform enjoy nice asymptotic properties we propose
the following procedure: first, make data continuous by adding a continuous
random noise and then apply the modified conditional distribution
transformation to get uniform iid series.

The first step can be called the \textit{continuous extension of a discrete
variable} which has been employed in different situations. For example
Ferguson (1967) uses some type of extension for simple hypothesis testing,
Denuit and Lambert (2005) and Neslezhova (2006) use it to apply a copulas
technique for discrete and discontinuous variables. The second step is the
probability integral transform (PIT) of the continued variables, which we will call \textit{randomized PIT}.
Resulting uniform iid series
can be tested using Bai (2003) or Corradi and Swanson (2006) tests. However,
in some cases these tests can not distinguish certain alternatives, so we
also propose test based on comparing joint empirical distribution functions
with the product of its theoretical uniform marginals by means of Cramer-von
Mises or Kolmogorov-Smirnov type statistics, developed by Kheifets (2011) for continuous distributions.

In a general setup, we do not know the true parameters, while the integral
transform using estimated parameters does not necessary provide iid uniform
random variates. Hence asymptotic properties and critical values of the
tests with estimated parameters have to be addressed. The estimation effect
changes the asymptotic distribution of the statistics and makes it \emph{%
data dependent}. Andrews (1997) proves that parametric bootstrap provides
correct critical values in this case using linear expansion of the
estimation effect, which arises naturally under the ML method. The idea of
orthogonal projecting the test statistics against the estimation effect due
to Wooldridge (1990) has been used in parametric moment tests, see Bontemps and
Meddahi (2006). The continuous version of the projection, often called
Khmaladze (1981) transformation, was employed in the tests of Koul and Stute
(1999) to specify the conditional mean, and of Bai and Ng (2001), Bai
(2003), Delgado and Stute (2008) to specify the conditional distribution.
These projection tests are not model invariant since they require to compute
conditional mean and variance derivatives, and also projections may cause a
loss in power. In this paper we apply a bootstrap  approach instead. In the case of ordered choice models an extensive
Monte Carlo comparison of specification tests has been done by Mora and Moro
(2008) in a static cross section context. They study two types of tests
based on moment conditions and on comparison of parametric and nonparametric
estimates.

Despite that there is some work on nonstationary discrete data models, cf.
Phillips and Park (2000), we stress stationary situations, but some ideas
could be extended to a more general set up as far as the conditional model
provides a full specification of the distribution of the dependent discrete
variable.

The contributions of this paper are following: 1) a new specification test
for dynamic discontinuous models is proposed, 2) we show that the test is invariant to the choice of distribution of the random noise added, 3) parameter estimation effect
of the test is studied, 4) under standard conditions we show the asymptotic properties of such tests, and 5) since asymptotic
distribution is case dependent, critical values can not be tabulated and we
prove that a bootstrap distribution approximation is valid.

The rest of the paper is organized as follows. Section 2 introduces
specification test statistics. Asymptotic properties
and bootstrap justification provided in Section 3. Monte Carlo experiments
are reported in Section 4.
Section~5 concludes.

\section{Test statistics}

In this section we introduce our goodness-of-fit statistics. Suppose that a
sequence of observations $(Y_{1},X_{1}),(Y_{2},X_{2}),...,(Y_{T},X_{T})$ is
given. Let $\Omega _{t}=\{X_{t},X_{t-1},\ldots ;Y_{t-1},Y_{t-2},\ldots \}$
be the information set at time $t$ (not including $Y_{t}$). We consider
a family of conditional cdf's $F(y|\Omega _{t},\theta )$, parameterized by $%
\theta \in \Theta $, where $\Theta \subseteq R^{L}$ is a finite dimensional
parameter space. We could allow for nonstationarity by permitting the change
in the functional form of the cdf of $Y_{t}$ using subscript $t$ in $F_{t}$.
Our null hypothesis of correct specification is \bigskip

$H_{0}$ : The conditional distribution of $Y_{t}$ conditional on $\Omega
_{t} $ is in the parametric family $F(y|\Omega _{t},\theta )$ for some $%
\theta _{0}\in \Theta $.

\bigskip

For example, for dynamic ordered discrete choice model the null hypothesis would mean that  $\exists \theta_0\in\Theta,\quad\forall j=0,\ldots,J,\quad \Pr\left(Y_t=j\right|\Omega_t)=p_j(\Omega_t,\theta_0),$ i.e. that all conditional probabilities are in a given parametric family.


For further analysis, we assume that the support of the conditional distributions $F(y|\Omega _{t},\theta )$ is a finite set of nonnegative integers $\{0,\ldots,J\}$ and $F(y|\Omega _{t},\theta )=\sum_{j\le y}\Pr_F(j|\Omega _{t},\theta )$, where $\Pr_F$ is the probability function at the discrete points.

The first step  is to obtain a \emph{%
continuous version} of $Y$.
 For any  random variable $Z\sim F_{z}$
with support in $[0,1]$ and $F_{z}$ continuous (but not necessary strictly increasing) define the \textit{continued by} $Z$
version of~$Y$,
\begin{equation*}
Y^{\dag }=Y+Z-1.
\end{equation*}%
Then the distribution of the continued version of $Y$ is
\begin{equation}\label{eq:Fstar}
F^{\dag }\left( y|\Omega
_{t}\right) =\Pr \left(Y^{\dag }\leq y|\Omega _{t}\right)=F\left([y]|\Omega
_{t}\right)+F_{z}(y-[y])\Pr \left([y]+1|\Omega _{t}\right),
\end{equation} which is strictly increasing on $%
\left[ -1,J\right] .$ The typical choice for $Z$ is the uniform in $\left[
0,1\right],$ so that
\begin{equation}\label{eq:FstarU}
F^{\dag }\left( y|\Omega _{t}\right) =F([y]|\Omega _{t})+(y-[y])\Pr
([y]+1|\Omega _{t}).
\end{equation}%
The binary choice case renders $F^{\dag }\left( y|\Omega _{t}\right)
=(y-[y])\left( 1-p_{t}\right) \ $for $y\in \lbrack -1,0)\ $and $F^{\dag
}\left( y|\Omega _{t}\right) =\left( 1-p_{t}\right) +(y-[y])p_{t}\ $for $%
y\in \left[ 0,1\right].$ Note, that $F^{\dag }$ coincides with $F$ in the domain of $F$.
We state next an "invariance property": for our purpose, it does not matter how to continue $Y$ and what distribution $F_z$ of the noise $Z$ to add. The unit support of $Z$ is needed to get a simple expression for $F^{\dag }$ in (\ref{eq:Fstar}), otherwise the resulting distribution will be a convolution
$F^{\dag}\left(y|\Omega _{t}\right)= \sum_{j=0}^J F_z\left(y+1-j\right)\Pr\left(j|\Omega_{t}\right)$. 
Continuation idea has been used to deal with discrete distributions, for example, to work with copulas with discrete marginals as in Denuit and
Lambert (2005).

The following proposition generalizes results about the probability integral transform. 
\begin{prop}
\label{propFY} (a) Under $H_{0}$
random variables $U_{t}=F^{\dag }(Y_{t}^{\dag }|\Omega _{t},\theta _{0})$
are iid uniform; (b) Invariant property of randomized PIT: realizations of $U_{t}$ are the same for any distribution $F_z$ in (\ref{eq:Fstar}) both under $H_{0}$ and under the alternative.
\end{prop}

Part (a) is a property of usual PIT with a continuous distribution $F^{\dag }$. Part (b) that realizations of 
$U_{t}$ are the same, means the following. Consider continuations of $Y_t$ by arbitrary $Z\sim F_z$ and uniform $Z_u\sim F_U$. Fix realizations $\{y_t\}$, $\{z_t\}$ and $\{z_{ut}\}$ from respective distributions. If $z_{ut}=F_z(z_t)$, then
\[
F^{\dag F_z }(y_{t}+z_t-1|\Omega _{t},\theta _{0}) = F^{\dag F_U }(y_{t}+z_{ut}-1|\Omega _{t},\theta _{0}),
\]
where $F^{\dag F_z }$ stresses dependence of $F^{\dag}$ on $F_z$ in (\ref{eq:Fstar}), $F^{\dag F_U }$ is as $F^{\dag}$ in (\ref{eq:FstarU}), continued by uniform, and $\Omega _{t}$ denotes here realized past.
Therefore, although a continued variable $Y_{t}^{\dag }$ and its distribution  $F^{\dag }$ depends on $F_z$, $F^{\dag }(Y_{t}^{\dag }|\Omega _{t},\theta _{0})$ is not and we can always use uniform variables $Z$ for continuation without affecting any properties of tests based on $U_t$.


Now we can use the fact that under the null hypothesis $U_{t}=F^{\dag
}(Y_{t}^{\dag }|\Omega _{t},\theta _{0}),$ $t=1,\ldots ,T$, are uniform on
[0,1] and iid random variables, so that $\Pr ({U}_{t-1}\leq r_{1},{%
U}_{t-2}\leq r_{2},...,{U}_{t-p}\leq r_{p})=r_{1}r_{2}\cdots r_{p}$, for $%
r=\left(r_1,\ldots,r_p\right)\in {[0,1]}^{p}$. This motivates us to consider the following empirical
processes
\begin{equation*}
{V}_{pT}(r)=\frac{1}{\sqrt{T-(p+1)}}\sum_{t=p+1}^{T}\left[ \prod_{j=1}^{p}I({%
U}_{t-j}\leq r_{j})-r_{1}r_{2}\ldots r_{p}\right].
\end{equation*}

If we do not know $\theta _{0}$ either $\{(Y_{t},X_{t}),t\leq 0\}$, we
approximate $U_{t}$ with $\hat{U}_{t}=F_{t}^{\dag }(Y_{t}^{\dag }|\tilde{%
\Omega}_{t},{\hat{\theta}})$ where ${\hat{\theta}}$ is an estimator of $%
\theta _{0}$ and the truncated information set is $\tilde{\Omega}%
_{t}=\{X_{t},X_{t-1},\ldots ,X_{1};Y_{t-1},Y_{t-2},\ldots ,Y_{1}\}$ and
write
\begin{equation}
\hat{V}_{pT}(r)=\frac{1}{\sqrt{T-(p+1)}}\sum_{t=p+1}^{T}\left[
\prod_{j=1}^{p}I(\hat{U}_{t-j}\leq r_{j})-r_{1}r_{2}\ldots r_{p}\right]
\label{eq:Vpn}
\end{equation}%
and
\begin{equation*}
D_{pT}=\Gamma (\hat{V}_{pT}(r))
\end{equation*}%
for any continuous functional $\Gamma (\cdot )$ from $\ell ^{\infty
}([0,1]^{p})$, the set of uniformly bounded real functions on $[0,1]^{p}$,
to $R$. In particular we use the Cramer-von Misses and Kolmogorov Smirnov test statistics
\begin{equation}
D_{pT}^{CvM}=\int_{[0,1]^{p}}\hat{V}_{pT}(r)^{2}dr\text{ or }%
D_{pT}^{KS}=\max_{[0,1]^{p}}\left\vert \hat{V}_{pT}(r)\right\vert .
\label{eq:Dpn}
\end{equation}
One further possibility is to test for $j$-lag
pairwise independence, using the process
\begin{equation}
\hat{V}_{2T,j}(r)=\frac{1}{\sqrt{T-j}}\sum_{t=j+1}^{T}\left[ I(\hat{U}%
_{t}\leq r_{1})I(\hat{U}_{t-j}\leq r_{2})-r_{1}r_{2}\right],
\label{eq:V2nj}
\end{equation}%
and corresponding test statistics $D_{2T,j}^{CvM}$ and $D_{2T,j}^{KS},$ say.

We can aggregate across $p$ or $j$ summing possibly with different weights $%
k(\cdot )$, obtaining generalized statistics
\begin{equation}
ADP_{T}=\sum_{p=1}^{T-1}k(p)D_{pT,}\text{ or }ADJ_{T}=%
\sum_{j=1}^{T-1}k(j)D_{2T,j}.  \label{eq:ADn}
\end{equation}

For $p=1$, $D_{1T}^{KS}$ delivers a generalization of Kolmogorov test to discrete distributions. Usually this test captures general deviations of marginal distribution but lacks power if only dynamics is misspecified. In particular, it does not have power
against alternatives where $U_{t}$ are
uniform on [0,1] but not independent.
For general $p$, $V_{pT}$ delivers a generalization of Kheifets (2011) to discrete distributions. This test should capture both deviations of marginal distribution and deviations in dynamics.

A more direct approach is based in Box-Pierce (1970) type of statistics, we
could consider%
\begin{equation*}
{BPU}_{m}:=
T\sum_{j=1}^{m}\hat{\rho}_{T,U}\left( j\right) ^{2},
\end{equation*}
$m=1,2,\ldots ,$ and $\hat{\rho}_{T,U}\left( j\right) $ are the sample
correlation coefficients of the $U_{t}^{\prime }s$ at lag $j.\ $Noting that $%
U_{t}$ should be uniform continuous iid random variables under the null of correctly
specified model, but might be correlated under alternative hypothesis of wrong
specification, ${BPU}_{m}$ is a good basis to design
goodness-of-fit tests. This idea is related to the tests of Hong (1998).
Alternatively, we can check autocorrelations of Gaussian residuals $\Phi(U_t)$
\begin{equation*}
{BPN}_{m}:=
T\sum_{j=1}^{m}\hat{\rho}_{T,\Phi(U)}\left( j\right) ^{2},
\end{equation*}
and normality of $\Phi(U_t)$ with Jarque-Bera test.
In addition we can check autocorrelations of discrete innovations,
\begin{equation*}
e_{t}=\frac{Y_{t}-E\left[ Y_{t}|\Omega _{t}\right] }{\left( \text{Var}\left[
Y_{t}|\Omega _{t}\right] \right) ^{1/2}},
\end{equation*}%
which are just the usual standardized probit residuals. We can define%
\begin{equation*}
{BPD}_{m}:=
T\sum_{j=1}^{m}\hat{\rho}_{T,e}\left( j\right) ^{2}
\end{equation*}%
and other statistics based on autocorrelations of squares of different types of  residuals.  The asymptotic distribution of these statistics
can be approximated by chi square distributions when the true parameters $%
\theta _{0}$ are known. Unlike tests based on empirical process, these tests can not capture some alternatives, for example if misspecification involves only higher order moments.

Parameter estimation affects the
asymptotic distribution of these statistics, as well as that of
those tests based on the empirical distribution of the $U_{t}^{\prime }s.$
There are different bootstrap and sampling techniques to approximate
asymptotic distribution, see for example Shao and Dongsheng (1995), Politis,
Romano and Wolf (1999). Since under $H_{0}$ we know the parametric conditional
distribution, we apply parametric bootstrap to mimic the $H_{0}$
distribution. We introduce the algorithm now for statistics  $\Gamma (\hat{V}_{2T})$.

\begin{enumerate}
\item Estimate model with initial data $\left( Y_{t},X_{t}\right) $, $%
t=1,2,...,T$, get parameter estimator $\hat{\theta}$, get test statistic $\Gamma (\hat{V}_{2T})$.

\item Simulate $Y_{t}^{\ast }$ with $F(\cdot |\Omega _{t}^{\ast },\hat{\theta%
})$ recursively for $t=1,2,...,T$, where the bootstrap information set is $%
\Omega _{t}^{\ast }=(X_{t},X_{t-1},\ldots ,Y_{t-1}^{\ast },Y_{t-2}^{\ast
},...)$.\label{basimulate}

\item Estimate model with simulated data $Y_{t}^{\ast }$, get $\theta ^{\ast
}$, get bootstrapped statistics $\Gamma (\hat{V}_{2T}^{\ast })$.

\item Repeat 2-3 $B$ times, compute the percentiles of the empirical
distribution of the $B$ boostrapped statistics.

\item Reject $H_{0}$ if $\Gamma (\hat{V}_{2T})$ is greater than the corresponding
$(1-\alpha )$th percentile.
\end{enumerate}

We will prove that $\Gamma (\hat{V}_{2T}^{\ast })$ has the same limiting
distribution as $\Gamma (\hat{V}_{2T})$. Bootstrapping other statistics is similar.

\section{Asymptotic properties of specification tests}

In this section we derive asymptotic properties of the statistics based on $V_{2T}$.  We start
with the simple case when we know parameters, then study how the
asymptotic distribution changes if we estimate parameters. We provide analyses under
the null, under the local and fixed alternatives. We first state all necessarily assumptions and propositions, then discuss them.

Let $\Vert \cdot \Vert $ denote
Euclidean norm for matrices, i.e. $\Vert A\Vert =\sqrt{\mathop{\rm tr}%
\nolimits(A^{\prime }A)}$ and for $\varepsilon >0,$ $B(a,\varepsilon )$\ is
an open ball in $R^{L}$ with the center in the point $a$ and the radius $%
\varepsilon $. In particular, for some $M>0$ denote $B_{T}=B\left( \theta
_{0},M T^{-1/2}\right) =\{\theta :||\theta -\theta _{0}||\leq M T^{-1/2}\}$.
For any discrete distributions $G$ and $F,$ with probability functions $\Pr_G$ and $\Pr_F$, and $r\in [0,1]$ define
\begin{eqnarray*}
d \left(G,F,r\right)&=&
G\left(F^{-1}\left(r\right)\right)-F\left(F^{-1}\left(r\right)\right)\\
&&+
\frac{r-F\left(F^{-1}(r)\right)}{\Pr_F\left(F^{-1}(r)+1\right)}\left(
\Pr_G\left(F^{-1}(r)+1\right)-\Pr_F\left(F^{-1}(r)+1\right)\right).
\end{eqnarray*}
We have $d \left(F,F,r\right)=0$, but $d \left(G,F,r\right)$ is not symmetric in $G$ and $F$.

\asss{1}{
Uniform boundedness away from zero: $\forall \varepsilon>0$, $\exists \delta>0$, such that $|F(0|\Omega _{t},\theta )|>\varepsilon$ and  $|F(j|\Omega _{t},\theta )-F(j-1|\Omega _{t},\theta )|>\varepsilon$ for $j=1,\ldots,J$ uniformly in $\theta\in B(\theta_0,\delta )$.
}

\asss{2}{Smoothness with respect to parameters:
\begin{itemize}
\item[(2.1)]
\begin{equation*}
E\max_{t=1,..,T}\sup_{u\in B_T}\max_{y}\left\vert
F\left(y |\Omega _{t}, u\right)
-F\left(y |\Omega _{t},\theta_0\right)\right\vert
=O\left( T^{-1/2}\right).
\end{equation*}
\item[(2.2)]  $\forall M\in(0,\infty)$, $\forall M_2\in(0,\infty)$ and $\forall\delta> 0 $
\begin{equation*}
\max_{y}
\frac{1}{\sqrt{T}}
\sum_{t=1}^{T}
\sup_{
\substack{
||u-v||\le
 M_2 T^{-1/2-\delta}\\
u,v\in B_T
}
 }
\left|
F\left(y |\Omega _{t}, u\right)
-F\left(y |\Omega _{t},v\right)
\right|
=o_{p}\left(1\right).
\end{equation*}
\item[(2.3)]  $\forall M\in(0,\infty)$,
there exists a uniformly continuous (vector) function $h(r)$
from $[0,1]^{2}$ to $R^{L}$, such that
\begin{equation*}
\sup_{v\in B_T}
\sup_{r\in \lbrack
0,1]^{2}}\left\vert
\frac{1}{\sqrt{T}}\sum_{t=2}^{T}h_t(r,v)-h(r)^{\prime }{\sqrt{T}\left(\theta_0-v\right) }%
\right\vert =o_{p}(1),
\end{equation*}
where
\begin{eqnarray*}
h_t(r,v)&=&
d \left(F\left(\cdot |\Omega _{t-1}, \theta_0\right),F\left(\cdot |\Omega _{t-1}, v\right),r_2\right) r_{1} \\
&&+
d \left(F\left(\cdot |\Omega _{t}, \theta_0\right),F\left(\cdot |\Omega _{t}, v\right),r_1\right)
I\left( F\left(Y_{t-1} |\Omega _{t-1}, \theta_0\right)\leq r_2\right).
\end{eqnarray*}
\end{itemize}
}


\asss{3}{Linear expansion of the estimator:
when the sample is generated by the null $F_{t}(y|\Omega
_{t},\theta _{0})$, the estimator $\hat{\theta}$ admits a linear expansion
\begin{equation}
\sqrt{T}(\hat{\theta}-\theta _{0})=\frac{1}{\sqrt{T}}\sum_{t=1}^{T}\ell
\left(Y_t,\Omega _{t}\right)+o_{p}(1),  \label{eqestlinear}
\end{equation}%
with $E_{F_{t}}\left( \ell
\left(Y_t,\Omega _{t}\right)|\Omega _{t}\right) =0$ and
$\frac{1}{T}\sum_{t=1}^{T}\ell
\left(Y_t,\Omega _{t}\right) \ell
\left(Y_t,\Omega _{t}\right)^{\prime
}\overset{p_{F_{t}}}{\rightarrow }\Psi .
$}

Dynamic probit/logit and general discrete choice models considered in Introduction can easily be adjusted to satisfy all  these assumptions.  Discrete support allows a simple analytical closed form of conditional distribution of continued variable by any continuous random variable on unit support as in (\ref{eq:FstarU}). Assumption 1 in particular requires that $F(0|\Omega _{t},\theta )$ and  $F(j|\Omega _{t},\theta )-F(j-1|\Omega _{t},\theta )$ for $j=1,\ldots,J$  are  bounded away from zero uniformly around $\theta_0$. To study parameter estimation effect we need to assume some smoothness of the distribution with respect to the parameter in Assumption 2 and a linear expansion of the estimator in Assumption 3. Note, the smoothness of the distribution with respect to the parameter is preserved after continuation, therefore Assumption 2 is similar to continuous case in Kheifets (2011); local Lipschitz continuity or existence of uniformly bounded first derivative of the distribution w.r.t. parameter is sufficient.
For bootstrap we will need  to strengthen Assumption 3 (see Assumption 3B below), although both conditions are standard
and satisfied for many estimators, for example for MLE. Note, that to establish  the convergence of the process $V_{2T}$ (with known $\theta_0$) under the null  (the following Proposition 2), we do not need these assumptions.

We now describe the asymptotic behavior of the process $V_{2T}({r})$ under
$H_{0}$. Denote by "$\weakc$" weak convergence of stochastic processes as random elements of
the Skorokhod space $D\left([0,1]^2\right)$.

\begin{prop}
\label{proplimitV} Under $H_0$
\begin{equation*}
V_{2T}\weakc V_{2\infty },
\end{equation*}%
where $V_{2\infty }({r})$ is bi-parameter zero mean Gaussian process with
covariance
\begin{equation*}
\mathop{\rm Cov}\nolimits_{V_{2\infty }}({r},{s})=(r_{1}\wedge
s_{1})(r_{2}\wedge s_{2})+(r_{1}\wedge s_{2})r_{2}s_{1}+(r_{2}\wedge
s_{1})r_{1}s_{2}-3r_{1}r_{2}s_{1}s_{2}.  
\end{equation*}
\end{prop}

To take into account the estimation effect on the asymptotic distribution, we use a Taylor expansion to approximate $\hat{V}_{2T}(%
{r})$ with $V_{2T}({r})$,
\begin{equation*}
\hat{V}_{2T}({r})=V_{2T}({r})+\sqrt{T}\left( {\hat{\theta}-\theta _{0}}%
\right) ^{\prime }h(r)+o_p(1)
\end{equation*}%
uniformly in $r$.
To identify the limit of $\hat{V}_{2T}(r)$, we need to
study limiting distribution of {$\sqrt{T}(\hat{\theta}-\theta _{0})$}, using the expansion from Assumption 3.
Define
\begin{equation*}
C_{T}(r,s,\theta )=E\left(
\begin{array}{c}
{V_{2T}(r)} \\
\frac{1}{\sqrt{T}}\sum_{t=1}^{T}\ell
\left(Y_t,\Omega _{t}\right)%
\end{array}%
\right) \left(
\begin{array}{c}
{V_{2T}(s)} \\
\frac{1}{\sqrt{T}}\sum_{t=1}^{T}\ell
\left(Y_t,\Omega _{t}\right)%
\end{array}%
\right) ^{\prime }
\end{equation*}%
and let $(V_{2\infty }(r),\psi _{\infty }^{\prime })^{\prime }$ be a zero
mean Gaussian process with covariance function $C(r,s,\theta
_{0})=\lim_{T\rightarrow \infty }C_{T}(r,s,\theta _{0})$. Dependence on $%
\theta $ on right hand side (rhs) comes through $U_{t}$ and $\ell\left(\cdot,\cdot\right)$.

Suppose the conditional distribution function $H(y|\Omega _{t})$ is not in
the parametric family $F(y|\Omega _{t},\theta )$ but has the same support.
For any $T_{0}\in
\{0,1,2,...,\}$ and $T\geq T_{0}$ define conditional on $\Omega _{t}$
conditional df%
\begin{equation*}
G_{T}(y|\Omega _{t},\theta )=\left( 1-\frac{\sqrt{T_{0}}}{\sqrt{T}}\right)
F(y|\Omega _{t},\theta )+\frac{\sqrt{T_{0}}}{\sqrt{T}}H(y|\Omega _{t}).
\end{equation*}%
Now we define local alternatives:

\bigskip

$H_{1T}$: Conditional cdf of $Y_{t}$ is equal to $G_{T}(y|\Omega _{t},\theta
_{0})$ with $T_{0}\neq 0$.

\bigskip Conditional cdf $G_{T}(y|\Omega _{t},\theta _{0})$
allow us to study all three cases: $H_{0}$ if $T_{0}=0$, $H_{1T}$ if $%
T=T_{0},T_{0}+1,T_{0}+2,...$ and $T_{0}\neq 0$ and $H_{1}$ if we fix $%
T=T_{0} $.
In the next proposition we provide the asymptotic distribution of
our statistics under the null and under the local alternatives.

\begin{prop}
\label{propall}

a) Suppose Assumptions 1-3 hold. Then under $%
H_{0}$
\begin{equation*}
\Gamma (\hat{V}_{2T})\overset{d}{\rightarrow }\Gamma (\hat{V}_{2\infty }),
\end{equation*}%
where $\hat{V}_{2\infty }({r})=V_{2\infty }(r)-h(r)^{\prime }\psi _{\infty
}. $

b) Suppose Assumptions 1-3 hold. Then under $%
H_{1T}$
\begin{equation*}
\Gamma (\hat{V}_{2T})\overset{d}{\rightarrow }\Gamma \left( \hat{V}_{2\infty
}+\sqrt{T_{0}}k-\sqrt{T_{0}}\xi ^{\prime }h\right) ,
\end{equation*}%
where
\begin{eqnarray*}
k(r) &=&\plim_{T\rightarrow \infty }\frac{1}{T}%
\sum_{t=2}^{T}\left\{
d \left(H\left(\cdot |\Omega _{t-1}\right),F\left(\cdot |\Omega _{t-1}, \theta_0\right),r_2\right) r_{1}
\right. \\
&&+\left.
d \left(H\left(\cdot |\Omega _{t}\right),F\left(\cdot |\Omega _{t}, \theta_0\right),r_1\right)
I\left( F\left(Y_{t-1} |\Omega _{t-1}, \theta_0\right)\leq r_2\right)
\right\} ,
\end{eqnarray*}%
and
\begin{equation}
\xi =\plim_{T\rightarrow \infty }\frac{1}{T}\sum_{t=1}^{T}\ell
\left(Y_t,\Omega _{t}\right).  \label{eqm}
\end{equation}

\end{prop}

Under $G_{T}$, the random variables $U_{t}=F^{\dag }(Y_{t}^{\dag }|\Omega _{t},\theta
_{0})$ are not anymore iid, instead $U_{t}^{\ast }=G_{T}^{\dag
}(Y_{t}^{\dag }|\Omega _{t},\theta _{0})$ are uniform iid.  The
first term in $k(r)$ controls for the lack of uniformity of $U_{t}$ (and it is similar to Bai's (2003) $k(r)$), it is zero when $U_{t}$
are uniform. The second term in $k(r)$ adds control for independence of $U_{t}$, cf. Kheifets (2011).

Under the alternative we may have also that (\ref{eqestlinear}) is not
centred around zero, since $E_{G_{T}}\left( \ell
\left(Y_t,\Omega _{t}\right)|\Omega _{t}\right) =\frac{\sqrt{%
T_{0}}}{\sqrt{T}}E_{H}\left( \ell
\left(Y_t,\Omega _{t}\right)|\Omega _{t}\right) $, therefore $\xi $
may be nonzero, which stands for information from estimation. This term does
not appear in Bai (2003) method, since his method projects out the estimation
effect. 

For the case of the one parameter empirical process, we can provide the
following corollary, which is similar to Bai (2003)'s single parameter
results.

\begin{cor}
a) Suppose Assumptions 1-3 hold. Then under $H_{0}$
\begin{equation*}
\Gamma (\hat{V}_{1T}\left( \cdot \right) )\overset{d}{\rightarrow }\Gamma (%
\hat{V}_{2\infty }(\cdot ,1)),
\end{equation*}%
where $\hat{V}_{1\infty }(\cdot )=V_{1\infty }(\cdot )-h(\cdot ,1)^{\prime
}\psi _{\infty }$ and $V_{1\infty }(\cdot )=V_{2\infty }(\cdot ,1).$

b) Suppose Assumptions 1-3 hold. Then under $H_{1T}$
\begin{equation*}
\Gamma (\hat{V}_{1T}\left( \cdot \right) )\overset{d}{\rightarrow }\Gamma (%
\hat{V}_{1\infty }\left( \cdot \right) +\sqrt{T_{0}}k_{1}\left(
\cdot \right) -\sqrt{T_{0}}h(\cdot ,1)^{\prime }\xi ),
\end{equation*}%
where for $r\in \left[ 0,1\right] $
\begin{equation*}
k_{1}(r)=\plim_{T\rightarrow \infty }\frac{1}{T}%
\sum_{t=2}^{T}
d \left(H\left(\cdot |\Omega _{t}\right),F\left(\cdot |\Omega _{t}, \theta_0\right),r\right).
\end{equation*}

\end{cor}
Note then that tests based on $\hat{V}_{1T}$ are not consistent against
alternatives for which $k_{1}=0$ and $h(\cdot ,1)=0$ but $k\neq 0$ or $h(\cdot,1)\neq 0$ on some set of positive measure.

We will justify our bootstrap procedure now, i.e. we prove that $\Gamma (\hat{V}_{2T}^{\ast })$ has the same limiting
distribution as $\Gamma (\hat{V}_{2T})$. We say that the sample is
distributed under $\{\theta _{T}:T\geq 1\}$ when there is a triangular array
of random variables $\{Y_{Tt}:T\geq 1,t\leq T\}$ with $(T,t)$ element generated by $%
F(\cdot |\Omega _{Tt},\theta _{T})$, where $\Omega
_{Tt}=(X_{t-1},X_{t-2},\ldots,Y_{Tt-1},Y_{Tt-2},\ldots)$. Similar arguments can be applied to other statistics.

\asss{3B} {
\label{asslb}For all nonrandom sequences $\{\theta _{T}:T\geq 1\}$ for
which $\theta _{T}\rightarrow \theta _{0}$, we have
\begin{equation*}
\sqrt{T}(\hat{\theta}-\theta _{T})=\frac{1}{\sqrt{T}}\sum_{t=1}^{T}
\ell
\left(Y_{Tt},\Omega _{Tt}\right)+o_{p}(1),
\end{equation*}%
under $\{\theta _{T}:T\geq 1\}$, where $E\left[ \ell
\left(Y_{Tt},\Omega _{Tt}\right)|\Omega
_{Tt}\right] =0$ and
\begin{equation*}
\frac{1}{T}\sum_{t=1}^{T}\ell
\left(Y_{Tt},\Omega _{Tt}\right) \ell
\left(Y_{Tt},\Omega _{Tt}\right)^{\prime
}\overset{p}{\rightarrow }\Psi .
\end{equation*}}

Note that the function $\ell
\left(\cdot,\cdot\right)$ now depends on $\theta_T$ and is assumed to be the same as in
Assumption 3. We require that estimators of close to $\theta _{0}$
points have \emph{the same} linear representation as the estimator of $%
\theta _{0}$ itself. 

\begin{prop}
\label{propGVhatB}Suppose Assumptions 1, 2 and 3B
hold. Then for any nonrandom sequence $\{\theta _{T}:T\geq 1\}$ for which $%
\theta _{T}\rightarrow \theta _{0}$, under $\{\theta _{T}:T\geq 1\}$,%
\begin{equation*}
\Gamma (\hat{V}_{2T}({r}))\overset{d}{\rightarrow }\Gamma (\hat{V}_{2\infty
}({r})).
\end{equation*}
\end{prop}

\section{Monte Carlo Simulation}

In this section we investigate the finite sample properties of our bootstrap tests
using Monte Carlo exercise. We use a simple dynamic Probit model
with one exogenous regressor with autoregressive dynamics. We consider
three specifications of dynamics
\begin{eqnarray*}
\ \text{Static model} &\text{:\ }&\pi _{t}=\pi _{0}+\beta X_{t}, \\
\text{Dynamic model}&\text{:\ }&\pi _{t}=\pi _{0}+\delta _{1}Y_{t-1}+\beta X_{t},\\
\text{Dynamic model with interactions} &:&\pi _{t}=\pi _{0}+\delta _{1}Y_{t-1}+\gamma
_{1}Y_{t-1}X_{t}+\beta X_{t},\ \  \gamma _{1}=-2\beta,
\end{eqnarray*}%
where in all specifications $X_{t}$ follows an AR$\left( 1\right) $ process,%
\begin{equation*}
X_{t}=\alpha _{1}X_{t-1}+e_{t},\ \ \ e_{t}\sim IIN\left( 0,1\right) ,
\end{equation*}%
and we set $\pi _{0}=0,\beta = 1,\delta _{1}=0.8,\alpha _{1}=0.8.$

We try 11 different scenarios of data generating processes (DGP) and null hypotheses (see Table \ref{t:scenarios}). In the first three we study the size properties of static, dynamic and dynamic with interactions probit models. Other scenarios check power when dynamics and/or marginals are misspecified. We take logit and $\left( \chi
_{1}^{2}-1\right) /2^{1/2}$ as alternative distributions.
We use sample sizes $T=100$ (Table \ref{t:T100}), $300$  (Table \ref{t:T300}) and $500$  (Table \ref{t:T500}) with 1000 replications. To estimate the
Bootstrap percentages of rejections we use a Warp bootstrap Monte Carlo (see
Giacomini, Politis and White, 2007) for all considered test statistics.  For tests based on "continued" residuals we
consider one-parameter ($p=1$) and two-parameter empirical processes ($p=2$) with $j=1$ and $j=2$ lags  and Cramer-von Misses (CvM) and Kolmogorov-Smirnov
(KS) criterions. To  make the results more readable, we denote them as $CvM_0=D_{1T}^{CvM}$, $CvM_1=D_{2T,1}^{CvM}$, $CvM_2=D_{2T,2}^{CvM}$ and $KS_0=D_{1T}^{KS}$, $KS_1=D_{2T,1}^{KS}$, $KS_2=D_{2T,2}^{KS}$.
We
consider Box-Pierce type tests for Gaussian and discrete residuals with $m=1,2,25.$ We also check normality of Gaussian residuals with a bootstrapped Jarque-Bera test~(JB). The results  of empirical process tests with further lags $j=3,4,5$ and correlation tests on uniform residuals do not provide additional information and are omitted.

\begin{table}
\caption{Different scenarios for Monte Carlo experiments.}
\label{t:scenarios}
\begin{center}
\setcounter{mycounter}{0}
\setlength{\tabcolsep}{1.8pt}\footnotesize\begin{tabular}{|c|ll|}
\hline
 & DGP & Null \tabularnewline
\hline
 \refstepcounter{mycounter} \label{c:psps}
\themycounter  & probit static & probit static\tabularnewline
 \stepcounter{mycounter}\themycounter  & probit dynamic & probit dynamic\tabularnewline
 \stepcounter{mycounter}\themycounter  & probit interactions & probit interactions\tabularnewline
\hline
\refstepcounter{mycounter} \themycounter\label{c:lsps}  & logit static & probit static\tabularnewline
 \refstepcounter{mycounter}\themycounter\label{c:csps}  & chi2 static & probit static\tabularnewline
 \refstepcounter{mycounter}\themycounter\label{c:ldps}    & logit dynamic & probit static\tabularnewline
 \refstepcounter{mycounter}\themycounter\label{c:cdps}    & chi2 dynamic & probit static\tabularnewline
 \refstepcounter{mycounter}\themycounter\label{c:lipd}    & logit interactions & probit dynamic\tabularnewline
 \refstepcounter{mycounter}\themycounter\label{c:cipd}   & chi2 interactions & probit dynamic\tabularnewline
 \refstepcounter{mycounter}\themycounter\label{c:lips}   & logit interactions & probit static\tabularnewline
 \refstepcounter{mycounter}\themycounter\label{c:cips}   & chi2 interactions & probit static\tabularnewline
\hline
\end{tabular}
\end{center}

\end{table}

\begin{table}
\caption{Percentage of rejections of test statistics with $T=100$.}
\label{t:T100}
\begin{center}
\setcounter{mycounter}{0}
\setlength{\tabcolsep}{1.8pt}\footnotesize\begin{tabular}{|c|r|ccccccccccccc|}
\hline
&  & $CvM_0 $  & $CvM_{1} $  & $CvM_{2} $  & $KS_0 $  & $KS_{1} $  & $KS_{2} $  & $BPN_{1} $  & $BPN_{2}$  & $BPN_{25} $  & $JB $  & $BPD_{1} $  & $BPD_{2} $  & $BPD_{25} $ \tabularnewline
\hline
   & $10\%$ &  $ 8.8 $  &  $ 7.4 $  &  $ 10.4 $  &  $ 8.4 $  &  $ 10.1 $  &  $ 9.2 $  &  $ 9.5 $  &  $ 9.6 $  &  $ 9.3 $  &  $ 8.3 $  &  $ 10.1 $  &  $ 10.6 $  &  $ 8.8 $ \tabularnewline
  \stepcounter{mycounter}\themycounter & $5\%$ &  $ 3.5 $  &  $ 4.3 $  &  $ 4.3 $  &  $ 3.9 $  &  $ 4.8 $  &  $ 4.7 $  &  $ 4.6 $  &  $ 3.7 $  &  $ 3.8 $  &  $ 4.4 $  &  $ 5.5 $  &  $ 5.1 $  &  $ 3.4 $ \tabularnewline
   & $1\%$ &  $ 0.3 $  &  $ 0.9 $  &  $ 0.4 $  &  $ 0.5 $  &  $ 1.1 $  &  $ 1.7 $  &  $ 1.5 $  &  $ 0.8 $  &  $ 0.3 $  &  $ 2.1 $  &  $ 0.8 $  &  $ 0.6 $  &  $ 0.3 $ \tabularnewline
\hline
   & $10\%$ &  $ 7.9 $  &  $ 8.3 $  &  $ 8.7 $  &  $ 7.0 $  &  $ 9.6 $  &  $ 9.2 $  &  $ 9.0 $  &  $ 10.6 $  &  $ 7.0 $  &  $ 11.2 $  &  $ 9.5 $  &  $ 10.7 $  &  $ 12.8 $ \tabularnewline
  \stepcounter{mycounter}\themycounter & $5\%$ &  $ 3.0 $  &  $ 3.6 $  &  $ 4.0 $  &  $ 2.8 $  &  $ 4.9 $  &  $ 4.5 $  &  $ 6.0 $  &  $ 4.0 $  &  $ 2.1 $  &  $ 5.9 $  &  $ 4.0 $  &  $ 4.4 $  &  $ 6.1 $ \tabularnewline
   & $1\%$ &  $ 0.0 $  &  $ 0.4 $  &  $ 0.1 $  &  $ 0.8 $  &  $ 1.2 $  &  $ 1.0 $  &  $ 0.9 $  &  $ 1.3 $  &  $ 0.3 $  &  $ 1.5 $  &  $ 0.4 $  &  $ 1.1 $  &  $ 0.7 $ \tabularnewline
\hline
   & $10\%$ &  $ 8.9 $  &  $ 10.0 $  &  $ 9.5 $  &  $ 7.7 $  &  $ 10.6 $  &  $ 9.4 $  &  $ 10.1 $  &  $ 11.3 $  &  $ 8.9 $  &  $ 10.7 $  &  $ 9.2 $  &  $ 9.2 $  &  $ 10.1 $ \tabularnewline
  \stepcounter{mycounter}\themycounter & $5\%$ &  $ 4.1 $  &  $ 4.1 $  &  $ 3.9 $  &  $ 3.6 $  &  $ 4.9 $  &  $ 5.0 $  &  $ 5.5 $  &  $ 5.5 $  &  $ 4.5 $  &  $ 5.4 $  &  $ 5.5 $  &  $ 3.8 $  &  $ 5.4 $ \tabularnewline
   & $1\%$ &  $ 0.1 $  &  $ 0.1 $  &  $ 0.2 $  &  $ 1.1 $  &  $ 0.5 $  &  $ 0.8 $  &  $ 1.2 $  &  $ 1.1 $  &  $ 0.5 $  &  $ 1.1 $  &  $ 0.6 $  &  $ 0.8 $  &  $ 0.5 $ \tabularnewline
\hline
\hline
   & $10\%$ &  $ 8.1 $  &  $ 9.0 $  &  $ 7.6 $  &  $ 8.4 $  &  $ 8.9 $  &  $ 9.9 $  &  $ 7.2 $  &  $ 8.8 $  &  $ 7.5 $  &  $ 9.9 $  &  $ 9.0 $  &  $ 9.2 $  &  $ 9.0 $ \tabularnewline
  \stepcounter{mycounter}\themycounter & $5\%$ &  $ 3.9 $  &  $ 4.6 $  &  $ 3.5 $  &  $ 3.6 $  &  $ 4.1 $  &  $ 3.7 $  &  $ 3.5 $  &  $ 4.1 $  &  $ 3.6 $  &  $ 3.0 $  &  $ 5.1 $  &  $ 4.6 $  &  $ 4.1 $ \tabularnewline
   & $1\%$ &  $ 0.5 $  &  $ 0.4 $  &  $ 0.3 $  &  $ 0.6 $  &  $ 0.6 $  &  $ 0.6 $  &  $ 1.2 $  &  $ 0.7 $  &  $ 0.6 $  &  $ 0.5 $  &  $ 1.0 $  &  $ 1.9 $  &  $ 0.7 $ \tabularnewline
\hline
   & $10\%$ &  $ 10.4 $  &  $ 9.5 $  &  $ 10.2 $  &  $ 12.0 $  &  $ 10.1 $  &  $ 11.1 $  &  $ 9.2 $  &  $ 11.5 $  &  $ 10.7 $  &  $ 20.3 $  &  $ 8.0 $  &  $ 7.5 $  &  $ 9.2 $ \tabularnewline
  \stepcounter{mycounter}\themycounter & $5\%$ &  $ 4.9 $  &  $ 6.1 $  &  $ 5.6 $  &  $ 5.9 $  &  $ 5.2 $  &  $ 5.7 $  &  $ 5.7 $  &  $ 6.3 $  &  $ 6.1 $  &  $ 12.6 $  &  $ 4.6 $  &  $ 3.7 $  &  $ 4.8 $ \tabularnewline
   & $1\%$ &  $ 0.5 $  &  $ 0.9 $  &  $ 0.3 $  &  $ 0.8 $  &  $ 1.2 $  &  $ 0.3 $  &  $ 1.0 $  &  $ 1.0 $  &  $ 0.7 $  &  $ 4.3 $  &  $ 1.8 $  &  $ 1.6 $  &  $ 0.9 $ \tabularnewline
\hline
   & $10\%$ &  $ 9.5 $  &  $ 11.0 $  &  $ 7.6 $  &  $ 9.2 $  &  $ 9.8 $  &  $ 9.3 $  &  $ 19.1 $  &  $ 15.4 $  &  $ 11.7 $  &  $ 11.0 $  &  $ 43.0 $  &  $ 35.3 $  &  $ 16.8 $ \tabularnewline
  \stepcounter{mycounter}\themycounter & $5\%$ &  $ 4.6 $  &  $ 4.9 $  &  $ 3.5 $  &  $ 3.5 $  &  $ 5.2 $  &  $ 4.6 $  &  $ 10.7 $  &  $ 9.0 $  &  $ 6.6 $  &  $ 4.7 $  &  $ 29.4 $  &  $ 20.5 $  &  $ 9.4 $ \tabularnewline
   & $1\%$ &  $ 0.4 $  &  $ 0.5 $  &  $ 0.8 $  &  $ 0.5 $  &  $ 1.4 $  &  $ 0.9 $  &  $ 2.9 $  &  $ 2.3 $  &  $ 0.8 $  &  $ 0.9 $  &  $ 11.0 $  &  $ 5.7 $  &  $ 2.9 $ \tabularnewline
\hline
   & $10\%$ &  $ 10.3 $  &  $ 10.9 $  &  $ 9.4 $  &  $ 9.2 $  &  $ 10.0 $  &  $ 9.3 $  &  $ 28.3 $  &  $ 26.4 $  &  $ 14.5 $  &  $ 13.7 $  &  $ 60.0 $  &  $ 50.6 $  &  $ 24.6 $ \tabularnewline
  \stepcounter{mycounter}\themycounter & $5\%$ &  $ 4.8 $  &  $ 5.2 $  &  $ 4.7 $  &  $ 3.9 $  &  $ 4.7 $  &  $ 5.2 $  &  $ 20.6 $  &  $ 16.4 $  &  $ 8.5 $  &  $ 7.7 $  &  $ 47.1 $  &  $ 37.0 $  &  $ 16.7 $ \tabularnewline
   & $1\%$ &  $ 0.1 $  &  $ 1.5 $  &  $ 0.1 $  &  $ 1.2 $  &  $ 1.3 $  &  $ 0.5 $  &  $ 9.4 $  &  $ 6.0 $  &  $ 2.6 $  &  $ 2.3 $  &  $ 26.0 $  &  $ 16.4 $  &  $ 5.6 $ \tabularnewline
\hline
   & $10\%$ &  $ 9.7 $  &  $ 9.2 $  &  $ 13.7 $  &  $ 9.9 $  &  $ 10.1 $  &  $ 13.0 $  &  $ 14.0 $  &  $ 26.6 $  &  $ 16.2 $  &  $ 9.9 $  &  $ 46.2 $  &  $ 57.4 $  &  $ 30.1 $ \tabularnewline
  \stepcounter{mycounter}\themycounter & $5\%$ &  $ 4.0 $  &  $ 3.7 $  &  $ 7.8 $  &  $ 3.6 $  &  $ 5.5 $  &  $ 7.8 $  &  $ 6.5 $  &  $ 18.8 $  &  $ 10.1 $  &  $ 3.8 $  &  $ 36.9 $  &  $ 45.1 $  &  $ 18.6 $ \tabularnewline
   & $1\%$ &  $ 0.8 $  &  $ 0.8 $  &  $ 2.5 $  &  $ 0.8 $  &  $ 1.1 $  &  $ 1.3 $  &  $ 0.6 $  &  $ 6.4 $  &  $ 2.3 $  &  $ 0.3 $  &  $ 17.1 $  &  $ 27.6 $  &  $ 5.0 $ \tabularnewline
\hline
   & $10\%$ &  $ 14.4 $  &  $ 16.9 $  &  $ 29.1 $  &  $ 16.0 $  &  $ 20.6 $  &  $ 34.4 $  &  $ 18.1 $  &  $ 55.7 $  &  $ 33.5 $  &  $ 20.0 $  &  $ 79.0 $  &  $ 82.2 $  &  $ 64.9 $ \tabularnewline
  \stepcounter{mycounter}\themycounter & $5\%$ &  $ 8.9 $  &  $ 10.0 $  &  $ 21.1 $  &  $ 9.9 $  &  $ 12.5 $  &  $ 18.5 $  &  $ 11.2 $  &  $ 48.4 $  &  $ 23.1 $  &  $ 12.4 $  &  $ 72.9 $  &  $ 81.0 $  &  $ 59.8 $ \tabularnewline
   & $1\%$ &  $ 0.9 $  &  $ 1.8 $  &  $ 3.9 $  &  $ 1.4 $  &  $ 4.2 $  &  $ 3.8 $  &  $ 3.4 $  &  $ 26.9 $  &  $ 11.6 $  &  $ 3.1 $  &  $ 58.1 $  &  $ 77.2 $  &  $ 43.8 $ \tabularnewline
\hline
   & $10\%$ &  $ 8.6 $  &  $ 14.7 $  &  $ 18.1 $  &  $ 7.6 $  &  $ 15.5 $  &  $ 13.0 $  &  $ 28.0 $  &  $ 42.0 $  &  $ 21.3 $  &  $ 9.2 $  &  $ 50.1 $  &  $ 79.8 $  &  $ 43.6 $ \tabularnewline
  \stepcounter{mycounter}\themycounter & $5\%$ &  $ 2.8 $  &  $ 8.5 $  &  $ 9.5 $  &  $ 3.8 $  &  $ 7.6 $  &  $ 6.6 $  &  $ 17.5 $  &  $ 29.0 $  &  $ 12.9 $  &  $ 3.9 $  &  $ 35.7 $  &  $ 69.8 $  &  $ 30.4 $ \tabularnewline
   & $1\%$ &  $ 0.6 $  &  $ 1.4 $  &  $ 2.1 $  &  $ 0.4 $  &  $ 1.3 $  &  $ 0.5 $  &  $ 5.3 $  &  $ 12.4 $  &  $ 4.4 $  &  $ 0.3 $  &  $ 22.4 $  &  $ 45.6 $  &  $ 10.7 $ \tabularnewline
\hline
   & $10\%$ &  $ 9.0 $  &  $ 28.1 $  &  $ 33.6 $  &  $ 8.1 $  &  $ 29.2 $  &  $ 28.4 $  &  $ 53.1 $  &  $ 85.1 $  &  $ 60.3 $  &  $ 8.8 $  &  $ 72.0 $  &  $ 99.9 $  &  $ 94.2 $ \tabularnewline
  \stepcounter{mycounter}\themycounter & $5\%$ &  $ 3.4 $  &  $ 17.6 $  &  $ 19.8 $  &  $ 3.3 $  &  $ 17.1 $  &  $ 11.8 $  &  $ 40.7 $  &  $ 76.3 $  &  $ 43.1 $  &  $ 5.3 $  &  $ 61.6 $  &  $ 99.7 $  &  $ 90.5 $ \tabularnewline
   & $1\%$ &  $ 0.2 $  &  $ 6.1 $  &  $ 4.2 $  &  $ 0.2 $  &  $ 3.4 $  &  $ 1.4 $  &  $ 23.2 $  &  $ 57.1 $  &  $ 22.0 $  &  $ 0.6 $  &  $ 40.8 $  &  $ 98.4 $  &  $ 72.3 $ \tabularnewline
\hline
\hline
\end{tabular}
\end{center}
\end{table}

\begin{table}
\caption{Percentage of rejections of test statistics with $T=300$.}
\label{t:T300}
\begin{center}
\setcounter{mycounter}{0}
\setlength{\tabcolsep}{1.8pt}\footnotesize\begin{tabular}{|c|r|ccccccccccccc|}
\hline
&  & $CvM_0 $  & $CvM_{1} $  & $CvM_{2} $  & $KS_0 $  & $KS_{1} $  & $KS_{2} $  & $BPN_{1} $  & $BPN_{2}$  & $BPN_{25} $  & $JB $  & $BPD_{1} $  & $BPD_{2} $  & $BPD_{25} $ \tabularnewline
\hline
   & $10\%$ &  $ 8.3 $  &  $ 9.2 $  &  $ 9.2 $  &  $ 9.1 $  &  $ 9.0 $  &  $ 10.4 $  &  $ 8.0 $  &  $ 8.3 $  &  $ 8.6 $  &  $ 8.0 $  &  $ 8.3 $  &  $ 8.2 $  &  $ 10.2 $ \tabularnewline
  \stepcounter{mycounter}\themycounter & $5\%$ &  $ 4.1 $  &  $ 4.7 $  &  $ 4.6 $  &  $ 5.2 $  &  $ 4.8 $  &  $ 4.6 $  &  $ 3.6 $  &  $ 3.9 $  &  $ 4.4 $  &  $ 2.7 $  &  $ 4.2 $  &  $ 3.1 $  &  $ 4.8 $ \tabularnewline
   & $1\%$ &  $ 0.7 $  &  $ 1.0 $  &  $ 0.6 $  &  $ 0.8 $  &  $ 0.6 $  &  $ 0.4 $  &  $ 0.6 $  &  $ 0.7 $  &  $ 1.1 $  &  $ 0.6 $  &  $ 0.6 $  &  $ 0.6 $  &  $ 1.0 $ \tabularnewline
\hline
   & $10\%$ &  $ 9.3 $  &  $ 8.9 $  &  $ 10.2 $  &  $ 9.3 $  &  $ 9.8 $  &  $ 11.6 $  &  $ 7.6 $  &  $ 10.0 $  &  $ 9.7 $  &  $ 9.6 $  &  $ 9.0 $  &  $ 6.9 $  &  $ 7.4 $ \tabularnewline
  \stepcounter{mycounter}\themycounter & $5\%$ &  $ 5.5 $  &  $ 4.4 $  &  $ 5.9 $  &  $ 5.2 $  &  $ 4.9 $  &  $ 6.4 $  &  $ 3.8 $  &  $ 4.6 $  &  $ 4.8 $  &  $ 3.8 $  &  $ 3.7 $  &  $ 3.6 $  &  $ 3.2 $ \tabularnewline
   & $1\%$ &  $ 1.0 $  &  $ 1.1 $  &  $ 0.9 $  &  $ 1.4 $  &  $ 1.1 $  &  $ 0.9 $  &  $ 0.5 $  &  $ 0.7 $  &  $ 1.4 $  &  $ 0.7 $  &  $ 0.5 $  &  $ 0.2 $  &  $ 0.6 $ \tabularnewline
\hline
   & $10\%$ &  $ 8.5 $  &  $ 12.3 $  &  $ 8.7 $  &  $ 8.7 $  &  $ 12.2 $  &  $ 9.9 $  &  $ 8.1 $  &  $ 9.9 $  &  $ 10.1 $  &  $ 9.0 $  &  $ 10.1 $  &  $ 9.4 $  &  $ 13.5 $ \tabularnewline
  \stepcounter{mycounter}\themycounter & $5\%$ &  $ 4.5 $  &  $ 5.1 $  &  $ 5.1 $  &  $ 3.3 $  &  $ 5.4 $  &  $ 5.3 $  &  $ 4.4 $  &  $ 5.2 $  &  $ 5.9 $  &  $ 4.2 $  &  $ 5.3 $  &  $ 4.2 $  &  $ 5.1 $ \tabularnewline
   & $1\%$ &  $ 1.1 $  &  $ 1.0 $  &  $ 1.5 $  &  $ 1.2 $  &  $ 0.6 $  &  $ 0.9 $  &  $ 0.9 $  &  $ 1.0 $  &  $ 1.3 $  &  $ 0.5 $  &  $ 0.6 $  &  $ 1.2 $  &  $ 1.5 $ \tabularnewline
\hline
\hline
   & $10\%$ &  $ 9.9 $  &  $ 10.2 $  &  $ 9.2 $  &  $ 9.5 $  &  $ 10.8 $  &  $ 9.8 $  &  $ 8.6 $  &  $ 9.3 $  &  $ 10.0 $  &  $ 11.2 $  &  $ 8.5 $  &  $ 9.6 $  &  $ 8.5 $ \tabularnewline
  \stepcounter{mycounter}\themycounter & $5\%$ &  $ 4.8 $  &  $ 4.9 $  &  $ 4.6 $  &  $ 5.3 $  &  $ 5.2 $  &  $ 4.9 $  &  $ 4.4 $  &  $ 5.1 $  &  $ 4.0 $  &  $ 6.1 $  &  $ 3.3 $  &  $ 3.2 $  &  $ 3.6 $ \tabularnewline
   & $1\%$ &  $ 0.9 $  &  $ 1.1 $  &  $ 0.6 $  &  $ 0.9 $  &  $ 1.3 $  &  $ 1.1 $  &  $ 0.5 $  &  $ 0.5 $  &  $ 1.0 $  &  $ 0.7 $  &  $ 0.8 $  &  $ 0.5 $  &  $ 0.9 $ \tabularnewline
\hline
   & $10\%$ &  $ 16.5 $  &  $ 15.1 $  &  $ 14.9 $  &  $ 15.8 $  &  $ 14.9 $  &  $ 14.4 $  &  $ 11.8 $  &  $ 11.2 $  &  $ 8.7 $  &  $ 41.4 $  &  $ 6.0 $  &  $ 7.2 $  &  $ 9.0 $ \tabularnewline
  \stepcounter{mycounter}\themycounter & $5\%$ &  $ 8.8 $  &  $ 7.9 $  &  $ 8.3 $  &  $ 9.6 $  &  $ 7.5 $  &  $ 8.5 $  &  $ 5.0 $  &  $ 6.1 $  &  $ 4.1 $  &  $ 30.4 $  &  $ 4.6 $  &  $ 6.3 $  &  $ 6.7 $ \tabularnewline
   & $1\%$ &  $ 1.6 $  &  $ 2.0 $  &  $ 1.8 $  &  $ 1.8 $  &  $ 2.2 $  &  $ 1.4 $  &  $ 0.9 $  &  $ 1.4 $  &  $ 1.0 $  &  $ 9.6 $  &  $ 3.0 $  &  $ 3.1 $  &  $ 2.9 $ \tabularnewline
\hline
   & $10\%$ &  $ 8.8 $  &  $ 15.4 $  &  $ 11.3 $  &  $ 9.0 $  &  $ 13.7 $  &  $ 10.1 $  &  $ 38.3 $  &  $ 29.9 $  &  $ 16.0 $  &  $ 9.6 $  &  $ 79.1 $  &  $ 69.0 $  &  $ 29.2 $ \tabularnewline
  \stepcounter{mycounter}\themycounter & $5\%$ &  $ 4.7 $  &  $ 9.8 $  &  $ 4.8 $  &  $ 5.4 $  &  $ 7.9 $  &  $ 6.1 $  &  $ 25.1 $  &  $ 21.8 $  &  $ 8.5 $  &  $ 5.8 $  &  $ 65.2 $  &  $ 58.1 $  &  $ 19.3 $ \tabularnewline
   & $1\%$ &  $ 0.5 $  &  $ 2.0 $  &  $ 0.3 $  &  $ 0.6 $  &  $ 1.5 $  &  $ 1.5 $  &  $ 11.8 $  &  $ 5.8 $  &  $ 1.5 $  &  $ 0.7 $  &  $ 43.7 $  &  $ 36.6 $  &  $ 6.8 $ \tabularnewline
\hline
   & $10\%$ &  $ 11.8 $  &  $ 12.2 $  &  $ 14.4 $  &  $ 12.3 $  &  $ 8.8 $  &  $ 13.2 $  &  $ 42.5 $  &  $ 35.1 $  &  $ 15.6 $  &  $ 24.5 $  &  $ 55.6 $  &  $ 47.1 $  &  $ 24.0 $ \tabularnewline
  \stepcounter{mycounter}\themycounter & $5\%$ &  $ 7.0 $  &  $ 5.9 $  &  $ 9.4 $  &  $ 6.2 $  &  $ 4.1 $  &  $ 6.3 $  &  $ 31.2 $  &  $ 25.0 $  &  $ 9.9 $  &  $ 16.3 $  &  $ 44.6 $  &  $ 39.5 $  &  $ 15.0 $ \tabularnewline
   & $1\%$ &  $ 1.1 $  &  $ 1.3 $  &  $ 1.9 $  &  $ 1.2 $  &  $ 1.2 $  &  $ 2.0 $  &  $ 15.3 $  &  $ 9.2 $  &  $ 2.6 $  &  $ 3.3 $  &  $ 35.3 $  &  $ 20.2 $  &  $ 5.8 $ \tabularnewline
\hline
   & $10\%$ &  $ 12.5 $  &  $ 18.3 $  &  $ 44.9 $  &  $ 13.0 $  &  $ 17.2 $  &  $ 44.0 $  &  $ 19.5 $  &  $ 67.0 $  &  $ 30.0 $  &  $ 12.6 $  &  $ 91.4 $  &  $ 96.5 $  &  $ 72.0 $ \tabularnewline
  \stepcounter{mycounter}\themycounter & $5\%$ &  $ 6.3 $  &  $ 12.5 $  &  $ 31.6 $  &  $ 6.4 $  &  $ 10.3 $  &  $ 28.9 $  &  $ 10.6 $  &  $ 55.7 $  &  $ 19.3 $  &  $ 6.2 $  &  $ 85.6 $  &  $ 94.0 $  &  $ 59.5 $ \tabularnewline
   & $1\%$ &  $ 1.0 $  &  $ 2.6 $  &  $ 9.8 $  &  $ 0.9 $  &  $ 2.7 $  &  $ 8.5 $  &  $ 2.7 $  &  $ 31.3 $  &  $ 7.4 $  &  $ 1.4 $  &  $ 71.4 $  &  $ 87.7 $  &  $ 37.3 $ \tabularnewline
\hline
   & $10\%$ &  $ 32.9 $  &  $ 42.5 $  &  $ 81.0 $  &  $ 29.5 $  &  $ 46.9 $  &  $ 88.8 $  &  $ 34.3 $  &  $ 92.0 $  &  $ 71.3 $  &  $ 24.8 $  &  $ 99.0 $  &  $ 99.7 $  &  $ 98.8 $ \tabularnewline
  \stepcounter{mycounter}\themycounter & $5\%$ &  $ 17.3 $  &  $ 25.9 $  &  $ 65.2 $  &  $ 19.7 $  &  $ 36.3 $  &  $ 80.2 $  &  $ 25.2 $  &  $ 88.8 $  &  $ 64.7 $  &  $ 17.6 $  &  $ 98.0 $  &  $ 99.6 $  &  $ 97.8 $ \tabularnewline
   & $1\%$ &  $ 3.0 $  &  $ 7.8 $  &  $ 36.3 $  &  $ 3.9 $  &  $ 14.3 $  &  $ 56.4 $  &  $ 12.3 $  &  $ 77.1 $  &  $ 42.1 $  &  $ 4.0 $  &  $ 94.5 $  &  $ 99.1 $  &  $ 95.6 $ \tabularnewline
\hline
   & $10\%$ &  $ 8.7 $  &  $ 33.1 $  &  $ 44.9 $  &  $ 9.7 $  &  $ 28.3 $  &  $ 33.7 $  &  $ 51.9 $  &  $ 83.2 $  &  $ 49.4 $  &  $ 9.3 $  &  $ 81.7 $  &  $ 99.6 $  &  $ 84.6 $ \tabularnewline
  \stepcounter{mycounter}\themycounter & $5\%$ &  $ 4.4 $  &  $ 22.4 $  &  $ 30.0 $  &  $ 4.7 $  &  $ 17.7 $  &  $ 21.4 $  &  $ 36.8 $  &  $ 74.5 $  &  $ 35.2 $  &  $ 5.5 $  &  $ 69.5 $  &  $ 98.8 $  &  $ 78.7 $ \tabularnewline
   & $1\%$ &  $ 1.1 $  &  $ 9.0 $  &  $ 11.1 $  &  $ 0.8 $  &  $ 6.8 $  &  $ 5.2 $  &  $ 18.2 $  &  $ 54.0 $  &  $ 17.0 $  &  $ 1.2 $  &  $ 37.9 $  &  $ 97.1 $  &  $ 48.1 $ \tabularnewline
\hline
   & $10\%$ &  $ 8.6 $  &  $ 46.2 $  &  $ 76.7 $  &  $ 9.1 $  &  $ 46.7 $  &  $ 76.9 $  &  $ 63.8 $  &  $ 99.5 $  &  $ 89.7 $  &  $ 11.0 $  &  $ 81.5 $  &  $ 100.0 $  &  $ 100.0 $ \tabularnewline
  \stepcounter{mycounter}\themycounter & $5\%$ &  $ 4.3 $  &  $ 32.8 $  &  $ 63.9 $  &  $ 4.3 $  &  $ 36.1 $  &  $ 67.8 $  &  $ 51.0 $  &  $ 98.8 $  &  $ 81.8 $  &  $ 5.6 $  &  $ 68.2 $  &  $ 100.0 $  &  $ 100.0 $ \tabularnewline
   & $1\%$ &  $ 0.4 $  &  $ 11.7 $  &  $ 37.3 $  &  $ 0.4 $  &  $ 18.4 $  &  $ 42.8 $  &  $ 31.7 $  &  $ 94.9 $  &  $ 63.9 $  &  $ 0.9 $  &  $ 39.4 $  &  $ 100.0 $  &  $ 99.3 $ \tabularnewline
\hline
\hline
\end{tabular}\end{center}
\end{table}

\begin{table}
\caption{Percentage of rejections of test statistics with $T=500$.}
\label{t:T500}
\begin{center}
\setcounter{mycounter}{0}
\setlength{\tabcolsep}{1.8pt}\footnotesize\begin{tabular}{|c|r|ccccccccccccc|}
\hline
&  & $CvM_0 $  & $CvM_{1} $  & $CvM_{2} $  & $KS_0 $  & $KS_{1} $  & $KS_{2} $  & $BPN_{1} $  & $BPN_{2}$  & $BPN_{25} $  & $JB $  & $BPD_{1} $  & $BPD_{2} $  & $BPD_{25} $ \tabularnewline
\hline
   & $10\%$ &  $ 10.2 $  &  $ 8.1 $  &  $ 9.2 $  &  $ 9.8 $  &  $ 8.6 $  &  $ 8.0 $  &  $ 11.9 $  &  $ 11.7 $  &  $ 11.1 $  &  $ 11.3 $  &  $ 10.1 $  &  $ 8.6 $  &  $ 9.8 $ \tabularnewline
  \stepcounter{mycounter}\themycounter & $5\%$ &  $ 4.7 $  &  $ 4.5 $  &  $ 4.9 $  &  $ 4.2 $  &  $ 3.9 $  &  $ 4.3 $  &  $ 6.0 $  &  $ 5.4 $  &  $ 5.6 $  &  $ 5.3 $  &  $ 5.2 $  &  $ 4.7 $  &  $ 4.7 $ \tabularnewline
   & $1\%$ &  $ 0.5 $  &  $ 1.0 $  &  $ 0.7 $  &  $ 0.7 $  &  $ 0.6 $  &  $ 1.0 $  &  $ 0.6 $  &  $ 0.8 $  &  $ 0.8 $  &  $ 1.2 $  &  $ 1.1 $  &  $ 0.4 $  &  $ 1.1 $ \tabularnewline
\hline
   & $10\%$ &  $ 9.1 $  &  $ 8.0 $  &  $ 8.3 $  &  $ 10.2 $  &  $ 8.6 $  &  $ 10.7 $  &  $ 11.6 $  &  $ 11.5 $  &  $ 8.4 $  &  $ 9.2 $  &  $ 8.6 $  &  $ 9.9 $  &  $ 10.6 $ \tabularnewline
  \stepcounter{mycounter}\themycounter & $5\%$ &  $ 4.3 $  &  $ 4.9 $  &  $ 4.4 $  &  $ 4.7 $  &  $ 3.8 $  &  $ 4.1 $  &  $ 5.6 $  &  $ 4.8 $  &  $ 3.9 $  &  $ 4.6 $  &  $ 4.3 $  &  $ 5.7 $  &  $ 5.6 $ \tabularnewline
   & $1\%$ &  $ 0.7 $  &  $ 0.8 $  &  $ 0.8 $  &  $ 0.5 $  &  $ 0.2 $  &  $ 0.6 $  &  $ 0.6 $  &  $ 1.4 $  &  $ 0.5 $  &  $ 0.5 $  &  $ 1.4 $  &  $ 0.8 $  &  $ 0.8 $ \tabularnewline
\hline
   & $10\%$ &  $ 9.1 $  &  $ 8.9 $  &  $ 9.9 $  &  $ 9.2 $  &  $ 8.6 $  &  $ 8.5 $  &  $ 10.0 $  &  $ 11.6 $  &  $ 10.7 $  &  $ 11.0 $  &  $ 10.7 $  &  $ 10.2 $  &  $ 10.0 $ \tabularnewline
  \stepcounter{mycounter}\themycounter & $5\%$ &  $ 4.1 $  &  $ 3.9 $  &  $ 3.1 $  &  $ 4.8 $  &  $ 4.5 $  &  $ 4.9 $  &  $ 5.7 $  &  $ 7.1 $  &  $ 5.4 $  &  $ 5.2 $  &  $ 4.6 $  &  $ 5.3 $  &  $ 4.8 $ \tabularnewline
   & $1\%$ &  $ 0.7 $  &  $ 1.1 $  &  $ 1.0 $  &  $ 1.0 $  &  $ 0.4 $  &  $ 0.4 $  &  $ 1.1 $  &  $ 1.0 $  &  $ 1.9 $  &  $ 1.5 $  &  $ 1.5 $  &  $ 2.0 $  &  $ 1.6 $ \tabularnewline
\hline
\hline
   & $10\%$ &  $ 10.7 $  &  $ 8.9 $  &  $ 10.9 $  &  $ 10.7 $  &  $ 10.8 $  &  $ 8.8 $  &  $ 11.5 $  &  $ 9.2 $  &  $ 10.2 $  &  $ 7.8 $  &  $ 10.0 $  &  $ 10.4 $  &  $ 12.3 $ \tabularnewline
  \stepcounter{mycounter}\themycounter & $5\%$ &  $ 5.2 $  &  $ 4.3 $  &  $ 4.8 $  &  $ 5.2 $  &  $ 3.9 $  &  $ 4.1 $  &  $ 6.1 $  &  $ 3.5 $  &  $ 5.2 $  &  $ 4.4 $  &  $ 5.1 $  &  $ 6.7 $  &  $ 7.2 $ \tabularnewline
   & $1\%$ &  $ 0.4 $  &  $ 0.8 $  &  $ 1.0 $  &  $ 0.7 $  &  $ 0.6 $  &  $ 0.4 $  &  $ 0.7 $  &  $ 1.1 $  &  $ 0.7 $  &  $ 0.7 $  &  $ 1.4 $  &  $ 1.3 $  &  $ 1.9 $ \tabularnewline
\hline
   & $10\%$ &  $ 17.1 $  &  $ 16.8 $  &  $ 17.0 $  &  $ 20.0 $  &  $ 19.4 $  &  $ 17.8 $  &  $ 11.8 $  &  $ 12.1 $  &  $ 8.5 $  &  $ 53.2 $  &  $ 7.9 $  &  $ 8.0 $  &  $ 14.0 $ \tabularnewline
  \stepcounter{mycounter}\themycounter & $5\%$ &  $ 11.4 $  &  $ 10.0 $  &  $ 9.7 $  &  $ 13.4 $  &  $ 12.0 $  &  $ 12.1 $  &  $ 4.4 $  &  $ 5.5 $  &  $ 3.8 $  &  $ 45.1 $  &  $ 5.6 $  &  $ 5.4 $  &  $ 9.1 $ \tabularnewline
   & $1\%$ &  $ 3.6 $  &  $ 3.7 $  &  $ 3.2 $  &  $ 4.5 $  &  $ 4.9 $  &  $ 3.3 $  &  $ 1.3 $  &  $ 1.7 $  &  $ 0.7 $  &  $ 23.6 $  &  $ 3.4 $  &  $ 3.7 $  &  $ 5.3 $ \tabularnewline
\hline
   & $10\%$ &  $ 8.7 $  &  $ 17.1 $  &  $ 9.8 $  &  $ 10.2 $  &  $ 15.5 $  &  $ 8.9 $  &  $ 46.1 $  &  $ 36.2 $  &  $ 17.6 $  &  $ 9.3 $  &  $ 88.7 $  &  $ 81.2 $  &  $ 42.0 $ \tabularnewline
  \stepcounter{mycounter}\themycounter & $5\%$ &  $ 5.3 $  &  $ 8.9 $  &  $ 4.9 $  &  $ 4.9 $  &  $ 6.9 $  &  $ 3.4 $  &  $ 32.9 $  &  $ 27.3 $  &  $ 11.8 $  &  $ 3.7 $  &  $ 82.9 $  &  $ 69.0 $  &  $ 29.5 $ \tabularnewline
   & $1\%$ &  $ 0.6 $  &  $ 1.5 $  &  $ 1.0 $  &  $ 0.9 $  &  $ 1.1 $  &  $ 0.8 $  &  $ 17.2 $  &  $ 8.9 $  &  $ 2.4 $  &  $ 0.6 $  &  $ 53.0 $  &  $ 46.8 $  &  $ 8.5 $ \tabularnewline
\hline
   & $10\%$ &  $ 15.2 $  &  $ 11.9 $  &  $ 14.8 $  &  $ 13.3 $  &  $ 11.3 $  &  $ 15.9 $  &  $ 39.9 $  &  $ 36.3 $  &  $ 18.0 $  &  $ 38.0 $  &  $ 53.7 $  &  $ 42.6 $  &  $ 21.5 $ \tabularnewline
  \stepcounter{mycounter}\themycounter & $5\%$ &  $ 8.1 $  &  $ 5.6 $  &  $ 10.8 $  &  $ 7.5 $  &  $ 4.8 $  &  $ 8.0 $  &  $ 28.5 $  &  $ 28.0 $  &  $ 10.6 $  &  $ 27.9 $  &  $ 41.0 $  &  $ 34.4 $  &  $ 16.5 $ \tabularnewline
   & $1\%$ &  $ 2.2 $  &  $ 1.3 $  &  $ 2.8 $  &  $ 2.5 $  &  $ 1.7 $  &  $ 3.0 $  &  $ 12.0 $  &  $ 9.4 $  &  $ 3.0 $  &  $ 8.2 $  &  $ 18.2 $  &  $ 19.5 $  &  $ 9.4 $ \tabularnewline
\hline
   & $10\%$ &  $ 22.6 $  &  $ 34.0 $  &  $ 90.1 $  &  $ 25.1 $  &  $ 35.9 $  &  $ 92.9 $  &  $ 23.7 $  &  $ 97.6 $  &  $ 76.7 $  &  $ 10.0 $  &  $ 99.9 $  &  $ 100.0 $  &  $ 99.6 $ \tabularnewline
  \stepcounter{mycounter}\themycounter & $5\%$ &  $ 12.6 $  &  $ 23.3 $  &  $ 83.8 $  &  $ 14.7 $  &  $ 25.9 $  &  $ 86.5 $  &  $ 16.7 $  &  $ 95.9 $  &  $ 65.1 $  &  $ 5.0 $  &  $ 99.7 $  &  $ 100.0 $  &  $ 99.2 $ \tabularnewline
   & $1\%$ &  $ 3.4 $  &  $ 7.6 $  &  $ 53.9 $  &  $ 4.0 $  &  $ 10.4 $  &  $ 65.0 $  &  $ 5.1 $  &  $ 87.6 $  &  $ 44.7 $  &  $ 1.1 $  &  $ 99.2 $  &  $ 99.9 $  &  $ 97.1 $ \tabularnewline
\hline
   & $10\%$ &  $ 56.1 $  &  $ 73.2 $  &  $ 98.8 $  &  $ 58.3 $  &  $ 78.1 $  &  $ 99.6 $  &  $ 62.6 $  &  $ 99.5 $  &  $ 96.9 $  &  $ 30.4 $  &  $ 100.0 $  &  $ 100.0 $  &  $ 100.0 $ \tabularnewline
  \stepcounter{mycounter}\themycounter & $5\%$ &  $ 39.9 $  &  $ 59.9 $  &  $ 97.0 $  &  $ 41.9 $  &  $ 69.9 $  &  $ 98.8 $  &  $ 51.2 $  &  $ 99.3 $  &  $ 95.4 $  &  $ 20.5 $  &  $ 100.0 $  &  $ 100.0 $  &  $ 100.0 $ \tabularnewline
   & $1\%$ &  $ 13.4 $  &  $ 38.7 $  &  $ 88.5 $  &  $ 16.1 $  &  $ 48.2 $  &  $ 95.6 $  &  $ 31.1 $  &  $ 98.6 $  &  $ 91.4 $  &  $ 10.7 $  &  $ 100.0 $  &  $ 100.0 $  &  $ 99.8 $ \tabularnewline
\hline
   & $10\%$ &  $ 9.9 $  &  $ 58.4 $  &  $ 90.5 $  &  $ 10.3 $  &  $ 52.8 $  &  $ 88.7 $  &  $ 74.8 $  &  $ 99.6 $  &  $ 93.6 $  &  $ 7.9 $  &  $ 98.7 $  &  $ 100.0 $  &  $ 100.0 $ \tabularnewline
  \stepcounter{mycounter}\themycounter & $5\%$ &  $ 5.2 $  &  $ 43.1 $  &  $ 82.2 $  &  $ 5.0 $  &  $ 38.2 $  &  $ 79.8 $  &  $ 65.2 $  &  $ 99.3 $  &  $ 89.4 $  &  $ 4.5 $  &  $ 95.9 $  &  $ 100.0 $  &  $ 100.0 $ \tabularnewline
   & $1\%$ &  $ 1.1 $  &  $ 14.3 $  &  $ 59.4 $  &  $ 0.6 $  &  $ 14.1 $  &  $ 39.6 $  &  $ 44.5 $  &  $ 97.3 $  &  $ 78.0 $  &  $ 1.3 $  &  $ 86.0 $  &  $ 100.0 $  &  $ 99.7 $ \tabularnewline
\hline
   & $10\%$ &  $ 10.0 $  &  $ 74.7 $  &  $ 99.1 $  &  $ 11.3 $  &  $ 81.3 $  &  $ 98.0 $  &  $ 86.8 $  &  $ 100.0 $  &  $ 100.0 $  &  $ 9.4 $  &  $ 99.6 $  &  $ 100.0 $  &  $ 100.0 $ \tabularnewline
  \stepcounter{mycounter}\themycounter & $5\%$ &  $ 3.6 $  &  $ 62.5 $  &  $ 96.7 $  &  $ 4.8 $  &  $ 71.3 $  &  $ 95.6 $  &  $ 78.9 $  &  $ 100.0 $  &  $ 99.9 $  &  $ 4.6 $  &  $ 98.6 $  &  $ 100.0 $  &  $ 100.0 $ \tabularnewline
   & $1\%$ &  $ 0.5 $  &  $ 38.5 $  &  $ 87.5 $  &  $ 0.8 $  &  $ 33.1 $  &  $ 80.7 $  &  $ 64.0 $  &  $ 100.0 $  &  $ 99.9 $  &  $ 1.3 $  &  $ 87.1 $  &  $ 100.0 $  &  $ 100.0 $ \tabularnewline
\hline
\hline
\end{tabular}\end{center}
\end{table}

Now we discuss the performance of empirical process based tests in comparison with traditional correlation tests. For $T=100$ almost all test statistics are slightly undersized (cases 1-3). The situation improves with larger $T$, and CvM statistics approach faster to nominal rates than KS. Overall, empirical size at $T=500$ is very good. The situation with power is not unambiguous. No test can capture static logit alternative to the null hypothesis of static probit model even  at $T=500$ (case \ref{c:lsps}). On the other hand, when static $\chi^2$ alternative to the null hypothesis of static probit is considered (case \ref{c:csps}), there is some power at $T=300$ which improves with $T=500$ for all empirical process based tests. Since under the null and under the alternative we have static models, correlation tests do not have power. Normality test (JB) is doing well only in the latter case.
 When there is a slight dynamic misspecification  added to logit (case \ref{c:ldps}), $CvM_1$ and $KS_1$ improve, but when it is added to $\chi^2$ our tests  and JB  doing worse (case \ref{c:cdps}). Correlation tests, on the contrary display power against these dynamic alternatives. When the alternative has dynamic interactions, and the null is a dynamic probit (cases \ref{c:lipd} and \ref{c:cipd}), all tests (but JB for logit) are doing well, and even better if  higher lags are taken into account. Finally, when dynamic interactions are taken versus static model (cases \ref{c:lips} and \ref{c:cips}), power is very good, and increases when more lags are considered. Exceptions are "marginal tests" $CvM_0$, $KS_0$ and JB. To summarize, dynamic misspecification can be captured well by empirical process statistics and correlation tests. Misspecification in marginals, on the contrary, can not be distinguished at all by correlation tests but  empirical process statistics, possibly multi-parameter, still
 work, although further research in improving power of these tests is needed.

To develop our omnibus type tests we introduce additional continuous noise. An important question is the effect of this noise on the power of the tests. Since correlation tests based on discrete residuals $BPD_j$ do not use additional noise, while correlation tests based on continuous residuals $BPN_j$ do, we can use the difference in rejection rates between these sets of statistics under dynamic misspecification as an indirect measure of the effect of the introduced noise, though correlation tests are not consistent against static alternatives. From our Monte Carlo simulations we see that for all scenarios we consider, correlation tests based on discrete residuals perform better, indicating that some power losses may indeed be attributed to the introduced noise.  To overcome this problem, we plan to develop tests for discrete models based on alternative transformations of the data without introducing additional noise, but still consistent against a wide range of nonparametric alternative hypotheses.

\section{Conclusion}

In this paper we have proposed new tests for checking goodness-of-fit of conditional
distributions in nonlinear discrete time series models. Specification of the
conditional distribution (but not only conditional moments) is important in
many macroeconomics and financial applications. Due to the parameter estimation
effect, the asymptotic distribution depends on the model and specific
parameter values. We show that our parametric bootstrap provides a good
approximation to asymptotic distributions and renders feasible and simple
tests. Monte Carlo experiments have shown that  tests based on empirical processes have power  if misspecification comes from  dynamics.   If misspecification affects  marginals alone, correlation tests are inconsistent, while tests based on empirical processes have some power. Comparing to the continuous case, we may conclude that there is a reduction of power due to the additional noise which distribution is known under the alternative too.

\pagebreak

\section*{Appendix}
\addcontentsline{toc}{section}{Appendix}
\proofs{Proposition \ref{propFY}} {
Part (a) is a property of dynamic PIT with a continuous conditional distribution $F_{t}^{\dag }$, the proof can be found in Bai (2003). Part~(b)
follows from the fact that (omitting dependence on $t$, $\Omega _{t}$ and $\theta$)
\begin{eqnarray*}
F^{\dag}\left(Y+Z-1\right)
&=&F\left([Y+Z-1]\right)
+Z^U\Pr\left([Y+Z]\right) \\
&=&F\left(Y-1\right)
+Z^U\Pr\left(Y\right), 
\end{eqnarray*}
where
\[
Z^U=F_z\left(Y+Z-1-[Y+Z-1]\right)=F_z\left(Z\right) 
\]
is uniform for any $Z\sim F_z$ continuous and with $[0,1]$ support, by the usual static PIT property.
Therefore, although a continued variable $Y^{\dag }$ and its distribution  $F^{\dag }$ depends on $F_z$, $F^{\dag }(Y^{\dag })$ does not.
}

\proofs{Propositions 2} {
Assumption 1 in Kheifets (2011) is satisfied automatically after applying continuation defined in (\ref{eq:FstarU}), therefore Proposition 1 of  Kheifets (2011) holds.
}

\proofs{Propositions 3} {Follows from Kheifets (2011), we need only to check  that Assumption 2 in Kheifets (2011) is satisfied.

 Let $r =F^{\dag}\left(y\right).$ Note that $[y]=F^{-1}(r)$ but $F\left([y]\right)=F\left(F^{-1}(r)\right)$ equals $r$ only when $y=[y]$. The inverse of  $F^{\dag}$ is
\begin{eqnarray*}
y &=&{\left(F^{\dag}\right)}^{-1}\left(r\right)
=[y]+\frac{r-F\left([y]\right)}{\Pr\left([y]+1\right)}
=[y]+1+\frac{r-F\left([y]+1\right)}{\Pr\left([y]+1\right)}\\
&=&F^{-1}(r)+\frac{r-F\left(F^{-1}(r)\right)}{\Pr\left(F^{-1}(r)+1\right)}.
\end{eqnarray*}
 Note also that $\left(r-F\left([y]\right)\right)/\Pr\left([y]+1\right)=y-[y]\in[0,1]$. Take distribution $G$ with the same support as $F$.  We have different useful ways to write $d \left(G,F,r\right)$:
\begin{eqnarray}
d \left(G,F,r\right) &=&\eta^{\dag} \left(r\right)-r=
G^{\dag}\left(\left(F^{\dag}\right)^{-1}\left(r\right)\right)-r
=G^{\dag}\left(y\right)-r \nonumber \\
&=&
G\left(\left[y\right]\right)-F\left(\left[y\right]\right)
+
\left(y-[y]\right)\left(
\Pr_G\left(\left[y\right]+1\right)-\Pr_F\left(\left[y\right]+1\right)\right)
\label{eq:dy}\\
&=&
G\left(\left[y\right]+1\right)-F\left(\left[y\right]+1\right) \nonumber\\
&&+
\left(y-[y]-1\right)\left(
\Pr_G\left(\left[y\right]+1\right)-\Pr_F\left(\left[y\right]+1\right)\right)
\label{eq:dy1}\\
&=&
G\left(F^{-1}\left(r\right)\right)-F\left(F^{-1}\left(r\right)\right) \nonumber\\
&&+
\frac{r-F\left(F^{-1}(r)\right)}{\Pr_F\left(F^{-1}(r)+1\right)}\left(
\Pr_G\left(F^{-1}(r)+1\right)-\Pr_F\left(F^{-1}(r)+1\right)\right). \label{eq:dr}
\end{eqnarray}
Thus, noting  that $\Pr\left(\cdot\right)$ is bounded away from zero, we have that Assumption 2 in this paper is sufficient for Assumption 2 in Kheifets (2011):
\begin{itemize}
\item[(K2.1)]
\begin{equation*}
E\sup_{t=1,..,T}\sup_{u\in B_T}\sup_{r\in [0,1]}\left\vert \eta^{\dag}  _{t}\left( r,u,\theta_0\right)-r \right\vert
=O\left( T^{-1/2}\right).
\end{equation*}
\item[(K2.2)]  $\forall M\in(0,\infty)$, $\forall M_2\in(0,\infty)$ and $\forall\delta> 0 $
\begin{equation*}
\sup_{r\in [0,1]}
\frac{1}{\sqrt{T}}
\sum_{t=1}^{T}
\sup_{
\substack{
||u-v||\le
 M_2 T^{-1/2-\delta}\\
u,v\in B_T
}
 }
\left|
\eta^{\dag}_{t}\left( r,u,\theta_0\right)
-
\eta^{\dag}_{t}\left(
r,v,\theta_0\right)
\right|
=o_{p}\left(1\right).
\end{equation*}
\item[(K2.3)]  $\forall M\in(0,\infty)$, $\forall M_2\in(0,\infty)$ and $\forall\delta> 0 $
\begin{equation*}
\sup_{|r-s|\le M_2 T^{-1/2-\delta} }
\frac{1}{\sqrt{T}}
\sum_{t=1}^{T}
\sup_{u\in B_T}
\left|
\eta^{\dag}
_{t}\left(r, u,\theta_0\right)
-
\eta^{\dag}_{t}\left(s, u,\theta_0\right)
\right|
=o_{p}\left(1\right).
\end{equation*}
\item[(K2.4)]  $\forall M\in(0,\infty)$,
there exists a uniformly continuous (vector) function $h(r)$
from $[0,1]^{2}$ to $R^{L}$, such that
\begin{equation*}
\sup_{u\in B_T }\sup_{r\in \lbrack
0,1]^{2}}\left\vert
\frac{1}{\sqrt{T}}\sum_{t=2}^{T}h_t-h(r)^{\prime }{\sqrt{T}\left( u-\theta_0\right) }%
\right\vert =o_{p}(1).
\end{equation*}
where
\begin{equation*}
h_t=
\left(\eta^{\dag} _{t-1}\left( r_{2},u,\theta_0\right)
-r_{2}\right) r_{1}
+\left( \eta^{\dag} _{t}\left( r_{1},u,\theta_0\right) -r_{1}\right)
I\left( F^{\dag} _{t-1}\left(  Y^{\dag} _{t-1}|u\right)\leq r_{2}\right).
\end{equation*}
\end{itemize}
For Part a), take
$d \left(F\left(\cdot |\Omega _{t}, \theta_0\right),F\left(\cdot |\Omega _{t}, \hat\theta\right)\right)$. Then (K2.1), (K2.2), (K2.4) follow from  (2.1), (2.2) and (2.3) because of representation (\ref{eq:dr}). If we compare (\ref{eq:dy}) and (\ref{eq:dy1}) we see that $d(\cdot)$ is not only continuous in $r$, but piece-wise linear, so (K2.3) is satisfied automatically.

For Part b), take
$d \left(G_T\left(\cdot |\Omega _{t}, \theta_0\right),F\left(\cdot |\Omega _{t}, \hat\theta\right)\right)$ and use the additivity of $d(\cdot)$ in the first arguments:
\begin{eqnarray*}
d \left(G_T\left(\cdot |\Omega _{t}, \theta_0\right),F\left(\cdot |\Omega _{t}, \hat\theta\right)\right)&=&
\left(1-\frac{\sqrt{T_0}}{\sqrt{T}}\right) d \left(F\left(\cdot |\Omega _{t}, \theta_0\right),F\left(\cdot |\Omega _{t}, \hat\theta\right)\right)\\
&&+
\frac{\sqrt{T_0}}{\sqrt{T}} d \left(H\left(\cdot |\Omega _{t}\right),F\left(\cdot |\Omega _{t}, \hat\theta\right)\right).
\end{eqnarray*}
}

\proofs{Propositions 5} {The proof is similar if we consider
$d \left(F\left(\cdot |\Omega _{t}, \theta_T\right),F\left(\cdot |\Omega _{t}, \hat\theta_T\right)\right)$ under $\{\theta_T:T\ge 1\}$.}

\newpage

\pagebreak

\end{document}